\documentclass[12pt]{amsart}
\usepackage{Preamble}

\title{Quantum modules of Semipositive Toric Varieties}
\author{Jae Hwang Lee}

\begin{document}

\begin{abstract}
    A smooth projective toric variety $X=X_\Sigma$ has a geometric quotient description $V \git T$. Using $2|1$-pointed quasimap invariants, one can define a quantum $H^*(T)$-module $QM(X)$, which deforms a natural module structure given by the Kirwan map $H^*(T) \rightarrow H^*(X)$. The Batyrev ring of $X$, defined from combinatorial data of the fan $\Sigma$, has its natural module structure given by the quotient of a polynomial ring, say BatM$(X)$. In this paper, we prove that $QM(X)$ and BatM$(X)$ are naturally isomorphic when $X$ is semipositive.
\end{abstract}
\maketitle
\tableofcontents
\setcounter{tocdepth}{1} 
\section{Introduction}
\subsection{Main result}  Let $X$ be a smooth projective toric variety with its geometric quotient description $V \git T$, where $V$ is an $n$-dimensional complex vector space, $T$ is an $r$-dimensional complex torus. In \cite{jquantum}, it was defined the \textit{quantum $H^*(T)$-module} on $H^*(X)$ via quasimap invariants with one light point: for $\sigma \in H^*(T),~\alpha \in H^*(X)$ and $\{T_i\},\{T^j\}$ a basis and its dual basis of $H^*(X)$ respectively,
\[\sigma \star \alpha : = \sum_{\beta,i}q^{\beta} T^i \langle \alpha, T_i \mid \sigma \rangle_{0,2|1,\beta}^{0+}, \]
where $H^*(T)$ denotes the cohomology of the stack quotient $[V/T]$. In the same paper, it was conjectured \cite[Conjecture 1.1]{jquantum} that the quantum module of $X$ is the same as the Batyrev module given by the quotient map from a polynomial ring to the Batyrev ring, when $X$ is \textit{semipositive}, i.e., the anticanonical divisor is numerically effective. In this paper, we prove this conjecture.
\begin{Thm}
\label{conj}
    For a smooth projective semipositive toric variety $V\git T$, the quantum $H^*(T)$-module structure of $V\git T$ coincides with a natural module structure of the Batyrev ring of $V\git T$.
\end{Thm}

\subsection{Quantum cohomology} Quantum cohomology is an interesting and useful object in both physics and mathematics. It implicitly appeared in \cite{lerche} by physicists in the name of chiral rings arising from an $N=2$ superconformal model by deforming the product structure of the ordinary cohomology ring. This deformation, coined as the \textit{quantum product}, was further investigated through topological quantum field theories introduced by Witten \cite{Witten1,Witten2} and defined in a physics paper \cite{vafa}. Mathematically rigorous definitions of the quantum product were given in \cite{ruan,mcduff} via symplectic geometry and in \cite{kontmanin} via algebraic geometry using Gromov--Witten invariants. Roughly speaking, Gromov--Witten invariants count the number of curves in a given space passing through some prescribed subspaces. These numbers are intriguing to both physicists and mathematicians due to their connection to both string theory and mirror symmetry. Remarkably, the associativity of the quantum product, which is equivalent to the Witten-Dijkgraaf-Verlinde-Verlinde (WDVV) equations, was used by Kontsevich to derive a recursive formula that counts the number of degree $d$ rational plane curves \cite{kontmanin}.

In general, the computation of quantum cohomology or quantum products is difficult since there is no universal way to compute Gromov--Witten invariants. Many methods have been introduced to compute the quantum product. The quantum cohomology of the projective space $\PP^n$ can be computed by using axioms in GW theory \cite[Example 8.1.2.1]{coxmirror}. For some projective bundles over $\PP^n$, Qin--Ruan computed its quantum cohomology using excess intersection theory \cite{ruanqin}. For a flag manifold, Givental--Kim \cite{bumsiggivtoda} and Kim \cite{bumsigtoda} computed its quantum cohomology and related it to a Toda lattice. Recently, Iritani--Koto used the quantum $D$-module approach to give a decomposition of the quantum cohomology of projective bundles \cite{iritaniprojbdl}. 

In particular, computing quantum cohomology can be done through Givental's Mirror theorem from \cite{givmirr} for complete intersections in a smooth projective toric variety $X$. There are two cohomological-valued functions: $J$-functions from Givental's graph space and $I$-functions from the GKZ system. These are relatively easy to compute compared to figuring out quantum products; $I$-functions are somewhat easier to compute than $J$-functions. There is a statement that generates a relation in quantum cohomology proved in \cite{givequiv,bumsighyper} (see also \cite[Thm 10.3.1]{coxmirror}): \textit{If a quantum differential operator annihilates the $J$-function of $X$, then there is a relation in the (small) quantum cohomology.} This Theorem is very useful when the quantum differential operator is explicitly known. However, knowing some relations in quantum cohomology does not give the entire quantum cohomology. Batyrev defined a ring associated to a smooth projective toric variety, say \textit{Batyrev (quantum) ring} and denote by $QH_{\text{Bat}}^{*}(X)$, by using combinatorial data of the fan \cite{batyrev}. The $I$-function of $X$ is annihilated by an operator that is defined from the GKZ system of $X$. This operator gives a relation in the Batyrev ring in the same way a quantum differential operator does in the quantum cohomology. As an application of Givental's mirror theorem, the $J$-function and the $I$-function are the same when $X$ is Fano \cite[Prop 11.2.17]{coxmirror}. It implies that the quantum cohomology ring $QH^*(X)$ is isomorphic to the Batyrev ring $QH_{\text{Bat}}^{*}(X)$, see \cite{giventalhomo} and \cite[Example 11.2.5.2]{coxmirror}.

It is no longer true that there is such an isomorphism between the quantum cohomology and the Batyrev ring without the Fano condition. Under a weaker condition, being semipositive, there are some counterexamples to such a phenomenon. In \cite{spielbergcount} and \cite[Example 11.2.5.2]{coxmirror}, they showed that the quantum cohomology and the Batyrev ring are not isomorphic for the Hirzebruch surface of type 2, write $\F_2$. Spielberg showed that the quantum cohomology and the Batyrev ring of the semipositive toric variety $\PP(\calO_{\PP^2}(3)\oplus \calO_{\PP^2})$ are not isomorphic \cite{spielbergop3}. Costa--Mir\'{o}-Roig generated an infinite family of toric projective vector bundles over $\PP^1$ where such failure occurs by showing that there is a quantum product with infinitely many quantum corrections \cite{costamiro2}. In contrast, we want to point out that the semipositive condition is not equivalent to the failure of having an isomorphism between quantum cohomology and the Batyrev ring. Spielberg verified that quantum cohomology and the Batyrev ring are the same for the semipositive toric bundle $\PP(\calO_{\PP^1}^{r-2}\oplus \calO_{\PP^1}(1) )$ with $r \geq 3$ \cite{spielbergop011}. Some studies have been done to clarify the relationship between the quantum cohomology and the Batyrev ring for the semipositive case in \cite{kwokmfd,iritaniseidel,guestD,mcduffseidel}. All these research include a change of coordinates, or a nontrivial mirror transform. One can ask a question: \textit{is there a natural enumerative interpretation of $QH_{\text{Bat}}^{*}(X)$ when a smooth projective toric variety $X$ is semipositive?}

\subsection{Quantum module} Inspired by the moduli of stable quotients \cite{MOP}, Ciocan-Fontanin--Kim defined and developed the theory of stable quasimaps \cite{toricquasi,gitquasi,quasimaporbi}. This theory is closely related to the theory of stable maps via the wall-crossing \cite{genus0,genushigh,zhou} along $\epsilon \in [0+,\infty]$: GW invariants at $\epsilon =\infty$ and quasimap invariants at $\epsilon =0+$. One can define a quantum deformation similar to the quantum product using 3-pointed quasimap invariants and compare the quasimap quantum cohomology ring with the Batyrev ring. Battistella observed this 3-pointed quasimap deformation and obtained a negative answer to whether this quasimap quantum cohomology of $\F_2$ is the same as the Batyrev ring of $\F_2$ \cite{lucathesis}; it was pointed out that a crucial reason for this phenomenon comes from the failure of the divisor equation, since the forgetful map does not define the universal curve over a moduli space of quasimaps.

In \cite{bigI}, they introduced light points as markings to the moduli spaces of quasimaps. The map forgetting a light point does define the universal curve, so that the divisor equation holds. Let $X$ be a smooth projective toric variety with its geometric quotient description $V \git T$, where $V$ is an $n$-dimensional complex vector space, $T$ is an $r$-dimensional complex torus. Recently, using the idea of light points, the author defined the quantum $H^*(T)$-module structure via quasimap invariants with one light point \cite{jquantum}, where $H^*(T)$ denotes the cohomology of the stack quotient $[V/T]$. The defining equations of this quantum module structure are similar to the associativity of the quantum product. The quantum $H^*(T)$-module deforms the module structure given by the Kirwan map $H^*(T) \rightarrow H^*(X)$ \cite{kirwan}. In the same paper, the author also computed the quantum $H^*((\C^*)^2)$-module structure of $\F_2$ applying the virtual localization \cite{GrabPand} to the quasimap moduli spaces. As a result, it was verified that the quantum $H^*((\C^*)^2)$-module of $\F_2$ is isomorphic to the Batyrev ring of $\F_2$ realized as a module given by the quotient of a polynomial ring. We want to point out that this isomorphism between the Batyrev module and the quantum module for $\F_2$ does not require any changes of coordinates. On the other hand, the quantum Kirwan map defined in \cite{woodward1} requires a change of coordinates (see \cite[Example 3.6]{woodwardminimal} for $\F_2$).

In this paper, we give a proof of Theorem \ref{conj}, which was conjectured in \cite{jquantum}, using the method from \cite{givequiv,bumsighyper}. This requires the Givental connection and the quantum differential equations arising from the divisor equation. As a consequence, one can show that the quantum $H^*(T)$-modules of \textit{all} smooth projective semipositive toric varieties are the same as the Batyrev modules. For example, there are finitely many such toric surfaces, see Appendix in \cite{kwoksurf}. Unfortunately, there is no smooth projective toric Calabi--Yau varieties. Our result provides an enumerative interpretation of the Batyrev modules for smooth projective toric semipositive varieties. This \textit{weakens} the Fano condition of the statement in \cite{giventalhomo} and \cite[Example 11.2.5.2]{coxmirror} that provides an isomorphism between the quantum cohomology and the Batyrev rings via the quantum module defined in \cite{jquantum}.

In section \ref{sec_pre}, we recall the Mori cone and the primitive relations for toric varieties. In section \ref{sec_quantum}, we give two properties of the descendant quasimap invariants with light points and recall the definition of the quantum $H^*(T)$-module for a toric variety. In section \ref{sec_qconn}, we prove the quantum differential equation arising from the divisor equation. In section \ref{sec_Iftn}, we verify a statement that generates a relation in the quantum $H^*(T)$-modules. In section \ref{sec_bat}, we prove Theorem \ref{conj}.

\textbf{Acknowledgements.} I would like to express my gratitude to Renzo Cavalieri and Mark Shoemaker, who have supported the completion of this project. This project was partially supported by Renzo Cavalieri’s NSF DMS 2100962, Mark Shoemaker’s NSF Grant 1708104. Additionally, I would like to acknowledge the Beijing International Center for Mathematical Research (BICMR) at Peking University for their generous and kind accommodation and resources that made this research possible.

\section{Mori cones of toric varieties}
\label{sec_pre}
We work over the complex numbers $\C$.
\subsection{Toric geometric quotients} Let $M$ be a lattice of rank $n$ and $N$ its dual lattice. Let $\Sigma \subset N_{\R}$ be a smooth complete fan and write $X_\Sigma$ for the associated smooth projective toric variety. Denote the set of 1-dimensional cones by $\Sigma(1)$. For a ray $\rho \in \Sigma(1)$, write $u_\rho$ for the minimal generator. There is an exact sequence
\begin{equation}
    \label{eqn_toricses}
    0 \longrightarrow  M \longrightarrow  \oplus_{\rho \in \Sigma(1)} \Z D_\rho \longrightarrow  \text{Pic}(X_\Sigma) \longrightarrow 0,
\end{equation}
where $m\in M$ maps to $\sum_\rho \langle m, u_\rho \rangle D_\rho$. From \eqref{eqn_toricses}, one can obtain a torus in $\C^{\Sigma(1)}$ given by
\[T=\{(t_\rho) \in (\C^*)^{\Sigma(1)} \mid \prod_\rho t_{\rho}^{\langle e_i,u_\rho \rangle}=1 \text{ for } 1 \leq i \leq n  \}.\]
The natural multiplication in $T$ induces $T$-action on $\C^{\Sigma(1)}$. A \textit{primitive collection} is a subset $P$ of rays in $\Sigma(1)$ satisfying (i) $P \not\subseteq \sigma(1)$ for all $\sigma\in \Sigma$; (ii) every proper subset of $P$ is contained in $\sigma(1)$ for some $\sigma \in \Sigma$. Define the irrelevant subset in $\C^{\Sigma(1)}$ by $Z_\Sigma:= \bigcup_{P} \mathbf{V}(x_\rho \mid \rho \in P)$ where $P$ varies all primitive collections. The toric variety $X_\Sigma$ has a toric geometric quotient description
\begin{equation}\
\label{eqn_geometricquo}
    X_\Sigma \simeq V \backslash Z_{\Sigma} \git T, 
\end{equation}
where $V:=\C^{\Sigma(1)}$. We will write such a quotient as $V \git T$.

\subsection{Mori cones} Let $X_\Sigma$ be a smooth projective toric variety. Let $N_1(X_\Sigma)$ be a finite dimensional $\R$-vector space obtained as the scalar extension of the set of numerically equivalent proper 1-cycles.  
One can see that $N_1(X_\Sigma) \simeq H_2(X_\Sigma,\R)$
and its dual is $\text{Pic} (X_\Sigma)_\R$. The Mori cone $\text{NE}(X_\Sigma) \subset N_1(X_\Sigma)$ is defined as the closure of the cone generate by classes of irreducible complete curves. Denote by $\text{NE}(X_\Sigma)_\Z$ for the the set of integral classes in the Mori cone.
We will drop $X_\Sigma$ in the notation $\text{NE}(X_\Sigma),~\text{NE}(X_\Sigma)_\Z$, and $N_1(X_\Sigma)$ if the context is clear.

\subsection{Primitive relations} There exists an exact sequence from \eqref{eqn_toricses}
\begin{equation}
    \label{eqn_toricses2}
    0 \longrightarrow  N_1 \longrightarrow  \R^{\Sigma(1)} \longrightarrow  N_\R \longrightarrow 0,
\end{equation}
where $[C]$ goes to $(D_\rho \cdot C)_{\rho}$ and $e_\rho \in \R^{\Sigma(1)}$ to $u_\rho$. Thus, $N_1$ can be interpreted as the set of linear relations among the minimal generators of $\Sigma(1)$. For a primitive collection $P=\{u_{\rho_1},\ldots, u_{\rho_k}\}$, the vector $\sum_{i=1}^{k}u_{\rho_i}$ lies in the relative interior of a minimal cone $\gamma \in \Sigma$. Then, there is a relation
\begin{equation}
    \label{eqn_primitiverelation}
    u_{\rho_1}+\cdots+ u_{\rho_k} - \sum_{\rho \in \gamma(1)}c_\rho u_\rho=0,
\end{equation}
for some $c_\rho \in \Q_{>0}$. This is called the \textit{primitive relation} of $P$. Since $\Sigma$ is smooth, $c_\rho$ are positive integers and $P \cap \gamma(1) = \varnothing$, see \cite[Def 2.8 and Prop 3.1]{batyreveff}. Therefore, we obtain the coefficient vector $\beta_P := (b_\rho)_{\rho}$ that represents an effective class in $H_2(X_\Sigma,\Z)$ given by
\begin{equation}
\label{eqn_primitivecoeff}
b_\rho = \begin{cases}
    1 & \rho \in P \\
    -c_\rho & \rho \in \gamma(1) \\
    0 & \rho \in \text{otherwise}.
\end{cases}
\end{equation}
Then, it satisfies $\sum_\rho b_\rho u_\rho =0$. We call it the \textit{primitive class} associated to the primitive collection $P$.  From \cite[Thm 2.15]{batyreveff}, the Mori cone of effective 1-cycles coincides with the cone generated by all primitive classes
\begin{equation}
\label{eqn_moriprimitive}
    \text{NE}(X_\Sigma) = \sum_{P:\text{primitive collection}} \R_+ \beta_P.
\end{equation}

\section{Quantum modules}
\label{sec_quantum}
\subsection{Stable toric quasimaps} Given a smooth projective toric variety $X_\Sigma$ and fix $g\geq 0,~m\geq 2,~ k\geq 0$ satisfying $2g-2+m\geq 0$, and $\beta \in H_2(X_\Sigma,\Z)$. It was defined the moduli space $Q_{g,m|k}(X_\Sigma, \beta)$ of stable toric quasimaps of degree $\beta$ with light points \cite{toricquasi,bigI}. Here, $m/k$ is the number of heavy/light points, respectively.
\begin{Def}
    We say $\beta \in H_2(X_\Sigma,\Z)$ is \textbf{(quasimap)-effective} if $\beta$ can be realized as a class of a $m|k$-stable genus $g$ stable toric quasimap to $X_\Sigma$. Denote the set of all (quasimap)-effective classes by $\text{Eff}(X_\Sigma)$ (or $\text{Eff}$ if the context is clear).
\end{Def}
The set $\text{Eff}(X_\Sigma)$ forms a semigroup with the property that if $\beta_1 + \beta_2 =0$ for $\beta_i \in $Eff$(X_\Sigma)$, then $\beta_1=\beta_2=0$ \cite[Lem 3.1.3]{toricquasi}.
\begin{Rem}
\label{remark}
\begin{enumerate}
    \item A quasimap consists of the data \[(C;x_1,\ldots, x_m;y_1,\ldots, y_k,\{L_\rho\}_{\rho \in \Sigma(1)},\{s_\rho\}_{\rho \in \Sigma(1)}),\] where $L_\rho$ are line bundles over a connected nodal curve $C$ and $s_\rho$ are sections of $L_\rho$. The line bundles $L_\rho$ can be trivialized through isomorphisms $\phi_l:\otimes_{\rho}L_{\rho}^{\otimes \langle l,\rho \rangle} \xrightarrow[]{\sim} \calO_\C$ for all $l \in M$.
    \item For a quasimap whose source curve is $C$, there exists a finite subset $B \subset C$, the set of \textit{base points}. Away from each base point, the sections $s_\rho$ of a quasimap define a map to $X_\Sigma$.
    \item There are no \textit{rational tails}, i.e., a nonzero degree component with up to one special point, due to the stability condition.
\end{enumerate}
\end{Rem}
We refer the readers to \cite{toricquasi,gitquasi,bigI,jquantum} for more extensive information on the quasimap moduli spaces.

The moduli space $Q_{g,m|k}(X_\Sigma, \beta)$ of quasimaps carries the perfect obstruction theory so that its virtual fundamental class exists. For the geometric quotient description $V \git T$ of $X_\Sigma$ in \eqref{eqn_geometricquo}, there is a corresponding stack quotient $[V/T]$. The cohomology of this stack quotient is given by
\[H^*([V/T]) = H^*([\text{pt}/T]) = \C[\sigma_1,\ldots, \sigma_r]. \]
We denote it by $H^*(T)$. There are well-defined evaluation maps for both heavy and light points.
\begin{align*}
    &\text{ev}_{i} : Q_{g,m|k}(X_\Sigma,\beta) \rightarrow X_\Sigma, \\
    &\hat{\text{ev}}_{j} : Q_{g,m|k}(X_\Sigma,\beta) \rightarrow [V/T].
\end{align*}
Denote by $\psi_i$ the psi-class from the $i$th cotangent line bundle at the $i$th heavy point. As in \cite{bigI}, the genus $g$ $m|k$-quasimap descendant invariants are defined as follows: for $\alpha_1 ,\ldots, \alpha_m  \in H^*(X_\Sigma)$ and $\sigma_1, \ldots, \sigma_k  \in H^*([V/T] )$,
\begin{equation}
\label{eqn_compatibility}
    \langle \alpha_1 \psi_{1}^{a_1}, \ldots, \alpha_m \psi_{m}^{a_m} \mid \sigma_1 \ldots, \sigma_k \rangle _{g,m|k,\beta}^{\text{quasi},X_\Sigma}:= \int_{[Q_{g,m|k}(X_\Sigma,\beta) ]^{\text{vir}}} \prod_{i=1}^{m} \text{ev}_{i}^{*}(\alpha_i) \psi_{i}^{a_i}   \prod_{j=1}^{k} \hat{\text{ev}}_{j}^{*}(\sigma_j).
\end{equation}
The superscript $^{\text{quasi},X_\Sigma}$ will be omitted when the context is clear.

\subsection{Divisor equation} Forgetting a light marking defines the universal curve over a moduli space of quasimaps. Thus, for the forgetful map $\pi_{k+1}:Q_{g,m|k+1}(X_\Sigma, \beta)\rightarrow Q_{g,m|k}(X_\Sigma, \beta)$, the compatibility of the virtual fundamental classes holds
\[[Q_{g,m|k+1}(X_\Sigma, \beta)] = \pi_{k+1}^*[Q_{g,m|k}(X_\Sigma, \beta)].\]
Denote by $\mathcal{L}_i$ the cotangent line bundle over $Q_{g,m|k+1}(X_\Sigma,\beta)$ whose fiber over a point $(C;\underline{x},\underline{y},\underline{L},\underline{s})$ is the cotangent space $T_{x_i}^{*}C$, and by $\mathcal{L}_i '$ the cotangent line bundle over $Q_{g,m|k}(X_\Sigma,\beta)$. Then, we have the following
\begin{equation}
\label{eqn_comparison1}
c_1(\mathcal{L}_i) = \pi_{k+1}^* c_1(\mathcal{L}'_i) + D(i|k+1,0;\text{rest},\beta),
\end{equation}
where $D(i|k+1,0;\text{rest},\beta)$ denotes the closure of the class of quasimaps whose source has components: $C_1$ containing $x_i$ and $y_{k+1}$ on which the degree restricted to this component is zero, and $C_2$ containing the rest of the markings. From \eqref{eqn_compatibility} and \eqref{eqn_comparison1}, one can derive the descendant version of the divisor equation: for $\Tilde{D} \in H^2([V/T])$ homogeneous,
\begin{align}
    \label{diveqn}
    \langle \alpha_1 \psi_{1}^{a_1}, \ldots, \alpha_m \psi_{m}^{a_m} \mid \sigma_1 \ldots, \sigma_k, \Tilde{D} &\rangle _{0,m|k+1,\beta}=\big( \int_{\beta}D \big) \langle \alpha_1 \psi_{1}^{a_1} \ldots, \alpha_m \psi_{m}^{a_m} \mid \sigma_1 \ldots, \sigma_k \rangle _{0,m|k,\beta} \\
    + &  
    \sum_{l=1}^{m}\langle \alpha_1 \psi_{1}^{a_1}, \ldots, \psi_{l}^{a_l - 1} D ,\ldots, \alpha_m \psi_{m}^{a_m} \mid \sigma_1 \ldots, \sigma_k \rangle _{0,m|k,\beta},
    \nonumber
\end{align}
where $D \in H^2(X_\Sigma)$ is the image of $\Tilde{D}$ via the Kirwan map. Indeed, there is an open locus $S$ in $\hat{\text{ev}}_{k+1}^{-1}(\Tilde{D})$ where the $k+1$th light marking does not meet any base points. The evaluation $\hat{\text{ev}}_{k+1}$ restricted to $S$ factors through $X_\Sigma$, so that the light insertion $\Tilde{D}$ becomes $D$ on $S$. This gives us that the forgetful map $\pi_{k+1}$ is generically of degree $\int_\beta D$, which is in the first term of the right-hand side of \eqref{diveqn}. For the second term, it is related to $D(i|k+1,0;\text{rest},\beta)$ in \eqref{eqn_comparison1}. Since the component $C_1$ is of degree zero, there is no base point on $C_1$. It follows that the evaluation map $\hat{\text{ev}}_{k+1}$ restricted to $D(i|k+1,0;\text{rest},\beta)$ maps to $X_\Sigma$, so that the light insertion $\Tilde{D}$ becomes $D$. Then, the restriction of the evaluation factors through the composition of the forgetful map $\pi_{k+1}$ restricted to $D(i|k+1,0;\text{rest},\beta)$ followed by $\text{ev}_{i}$, since a quasimap in $D(i|k+1,0;\text{rest},\beta)$ is constant on the component $C_1$. Therefore, we obtain the insertion $D$ in the second term of the right-hand side of \eqref{diveqn} from \eqref{eqn_comparison1}. The rest of the argument is similar to that in GW theory.

\subsection{Genus zero topological recursion} Let $m \geq 2$, $k\geq 0$. Denote by $\mathcal{L}_1$ and $\mathcal{L}_1 '$ the cotangent line bundles over $Q_{g,m|k}(X_\Sigma,\beta)$ and $\overline{M}_{g,m|k}$ with respect to the first heavy point, respectively. For the forgetful map $\pi:Q_{0,2|k}(X_\Sigma,\beta) \rightarrow \overline{M}_{0,2|1} =pt$, a quasimap version of the comparison lemma holds
\begin{equation}
    \label{comparisonlemma}
    \mathcal{L}_1 = \pi^* \mathcal{L}_1 ' + \sum_{\substack{K\sqcup L = \{1,\ldots,k\},~1\in K \\ \beta_1 + \beta_2 = \beta}}D(1|K,\beta_1 ; 2|L,\beta_2),
\end{equation}
where $D(1|\{1\}\cup K,\beta_1 ; 2|L,\beta_2)$ is the divisor on $Q_{0,2|k}(X_\Sigma,\beta)$ whose source curve has one component with the first heavy point, the first light point, and all other light points marked by $K$, and the other component with the second heavy point, and all other light points marked by $L$. Recall that $c_1(\mathcal{L}_1) = \psi_1$ and denote $c_1(\mathcal{L}_1 ')$ by $\psi_1 '$. Let $\{T_a\}$ be a basis for $H^*(X_\Sigma)$, $\{T^a\}$ the dual basis and $\sigma_1,\ldots,\sigma_k \in H^*([V/T])$. We multiply $\text{ev}_1^* (\alpha_1) \psi_{1}^{a_1} \text{ev}_2^* (\alpha_2) \psi_{2}^{a_2} \hat{\text{ev}}_1(\sigma)$ to \eqref{comparisonlemma}. Since $\phi_1' =0$, by the splitting axiom in \cite{jquantum}, the following holds:
\begin{align}
    \label{trr}
    & \langle \alpha_1 \psi_{1}^{a_1+1}, \alpha_2 \psi_{2}^{a_2} \mid \sigma_1,\ldots, \sigma_k \rangle_{0,2|k,\beta}  \\
    &= \sum_{\substack{K\sqcup L = \{1,\ldots,k\},1\in K \\ \beta_1 + \beta_2 = \beta}}\sum_{a}
    \langle T_a, \alpha_2 \psi_{2}^{a_2} \mid \sigma_1, \sigma_{k_2},\ldots, \sigma_{k_{|K|}} \rangle_{0,2||K|,\beta_1}
    \langle \alpha_1 \psi_{1}^{a_1}, T^a \mid \sigma_{l_1},\ldots, \sigma_{l_{|L|}} \rangle_{0,2||L|,\beta_2},
    \nonumber
\end{align}
where $K=\{k_i\}$ with $k_1=1$ and $L=\{l_j\}$.

For the divisor equation and the topological recursion, we refer the readers to \cite[\S10.1-2]{coxmirror} or \cite[\S26.3-4]{mirrsymm}.

\subsection{Quantum $H^*(T)$-module structure} Fix $g=0$. Let $\{T_i\}$ be a basis for $H^*(X_\Sigma)$ and $\{T^j\}$ be the dual basis. In \cite{jquantum}, for $\sigma \in H^*([V/T])=H^*(T) $ and $\alpha \in H^*(X_\Sigma)$, the following was defined
    \begin{equation}
    \label{def_staroper}
        \sigma \star \alpha:= \sum_{\beta}\sum_{i} q^{\beta} \langle \alpha, T_i \mid \sigma \rangle_{0,2|1,\beta} T^i,
    \end{equation}
    where $q$ is the Novikov variable, and $\beta$ varies in $\text{Eff}(X_\Sigma)$. This action defines the \textit{quantum $H^*(T)$-module} structure on $H^*(X_\Sigma)$ via the equation
\begin{equation}
    \label{eqn_WDVV}
    \tau \star (\sigma \star \alpha)  = (\tau \cup \sigma) \star \alpha, 
\end{equation}
where $\tau \in H^*([V/T])$ and $\cup$ is the cup product in $H^*([V/T])$. Denote by $QM(X_\Sigma)$ for the quantum $H^*(T)$-module of $X_{\Sigma}$. We will frequently say the \textit{quantum module} of $X_\Sigma$. From the equation \eqref{eqn_WDVV}, to have full description of the module structure, it is enough to know $\sigma_k \star T_i$ for a homogeneous basis $\{\sigma_k\}$ of $H^0(T) \oplus H^2 (T)$.

\section{Quasimap Givental connection}
\label{sec_qconn}
In this section, we define an analogue of the Givental connection based on our quasimap module structure and flat sections of this connection. Then, some useful lemmas will be provided.

\subsection{Givental connections} 
Let $\{ T_i \}_{i=0}^{m}$ be a homogeneous basis for $H^*(X_\Sigma)$ with $T_0:=[X_\Sigma] \in H^0(X_\Sigma)$ and $r$ the rank of $H^2(X_\Sigma)$. Assume that $T_i$ are in $H^2(X_\Sigma)$ for $i=1,\ldots,r$. Since our toric variety $X_\Sigma$ is smooth, $H^2(X_\Sigma,\Z)$ is free. Via the Kirwan map $H^*(T) \rightarrow H^*(X_\Sigma)$, one can identify $H^0(T) \oplus H^2(T)$ with $H^0(X_\Sigma) \oplus H^2(X_\Sigma)$. Thus, there are unique lifts $\Tilde{T_i} \in H^0(T) \oplus H^2(T)$ for $T_i$. Denote by $t_i$ for the coordinates with respect to $T_i$. We consider the Givental connection only on the restriction $M:=H^0(X_\Sigma) \oplus H^2(X_\Sigma)$.
\begin{Def}
    Define the \textbf{quasimap Givental connection} $\nabla$ as a formal connection on the trivial bundle $M \times H^*(X_\Sigma)$ over $M$ as follows:
    \[
    \nabla_{\frac{\partial}{\partial t_i}} \big(\sum_{j=0}^{m}a_jT_j\big):=\hbar \sum_{j=0}^{m}\frac{\partial a_j}{\partial t_i}T_j - \sum_{j=0}^{m}a_j\Tilde{T_i} \star T_j,
    \]
    where $\hbar$ is a formal variable, $a_j$ are $H^*(X,\C) \otimes_\C \C[[q^\beta \mid q \in \text{Eff}]]$-valued functions of $t_0, \ldots, t_r$, and $i=0,\ldots, r$. We write $\nabla_i$ for $\nabla_{\frac{\partial}{\partial t_i}}$.
\end{Def}
Thus, the definition of a formal section of the trivial bundle $M \times H^*(X_\Sigma)$ over $M$ is given as follows. Let $\delta:=\sum_{i=1}^{r}t_iT_i$. Write $q^{\beta}=e^{\int_{\beta} \delta } = \prod_{i}e^{t_{i}\int_{\beta}T_{i} }$. It is easy to see that $\frac{\partial}{\partial t_i} q^\beta = q^\beta\big(\int_{\beta} T_i\big)$.
\begin{Def}
    For $a=0,\ldots,m$, define a formal section
    \[
    s_a:=e^{t_0/\hbar} \bigg( e^{\delta/\hbar} \cup T_a + \sum_{0\neq \beta \in \text{Eff}}\sum_{b=0}^{m} q^{\beta}\bigg\langle \frac{e^{\delta/\hbar} \cup T_a}{\hbar - \psi},  T_b  \bigg\rangle_{0,2,\beta} T^b   \bigg) .
    \]
\end{Def}

There is a quantum differential equation corresponding to the quasimap Givental connection $\nabla$ and its solutions are given by the formal sections $s_a$.
\begin{Pro}[Quantum differential equation]
\label{pro_qde}
The following equation holds
    \begin{equation}
    \label{eqn_qde}
    \hbar\frac{\partial s_a}{\partial t_i}  = \Tilde{T_i}\star s_a, \quad i=0,\ldots,r.
    \end{equation}
Equivalently, the formal section $s_a$ defines a flat section of the quasimap Givental connection $\nabla$, i.e., $\nabla_i s_a=0$ for all $i,a$.
\end{Pro}
\begin{proof}
    It is enough to verify the equation \eqref{eqn_qde} when $t_0=0$ since $\hbar \frac{\partial}{\partial t_0} s_a = \Tilde{T_0} \star s_a = s_a$. This is a consequence of the existence of the forgetful map at a light point.
    
    We start from
    \begin{equation*}
        s_a= e^{\delta/\hbar} \cup T_a + \sum_{0\neq \beta}\sum_{b=0}^{m} q^{\beta}\bigg\langle \frac{e^{\delta/\hbar} \cup T_a}{\hbar - \psi},  T_b  \bigg\rangle_{0,2,\beta} T^b  .
    \end{equation*}
    By using the Taylor expansion to $e^{\delta / \hbar}$ in the second term and using the geometric series to $\frac{1}{\hbar -\psi} = \hbar^{-1}\frac{1}{1 -(\psi/\hbar)} $, one can write
    \begin{align}
    \label{pro_qde_secexpn}
        &s_a = e^{\delta/\hbar} \cup T_a \\
        &+\sum_{\beta \neq 0}q^\beta \sum_b T^b \hbar^{-1} \sum_{n\geq 0} \sum_{k+l=n} \frac{\hbar^{-l-k}}{k!} \sum_{l_1+\cdots + l_r=k} \binom{k}{l_1,\ldots,l_r}  \prod_{j=1}^{r}t_{j}^{l_j}\langle T_a \psi^l\prod_{j=1}^{r}T_{j}^{l_j}  , T_b\rangle_{0,2,\beta} . \nonumber
    \end{align}

    We compute $\hbar\frac{\partial s_a}{\partial t_i}$. For the first term on the right-hand side of \eqref{pro_qde_secexpn},
    \begin{equation}
        \label{pro_qde1} 
        \hbar\frac{\partial}{\partial t_i } e^{\delta/\hbar} \cup T_a = e^{\delta/\hbar} \cup T_a \cup T_i.
    \end{equation}
    We apply $\hbar\frac{\partial}{\partial t_i }$ to the second term on the right-hand side of \eqref{pro_qde_secexpn}. By the Leibnitz rule,
    \begin{align}
    \label{pro_qde2}
        &\sum_{\beta \neq 0} q^\beta \big(\int_{\beta}T_i \big)\sum_b T^b \sum_{n\geq 0}\sum_{k+l=n}\frac{\hbar^{-l-k}}{k!} \sum_{l_1+\cdots + l_r=k} \binom{k}{l_1,\ldots,l_r} \prod_{j=1}^{r}t_{j}^{l_j} \langle T_a \psi^l\prod_{j=1}^{r}T_{j}^{l_j}  , T_b\rangle_{0,2,\beta}  \\
        +&\sum_{\beta \neq 0} q^\beta \sum_b T^b \sum_{n\geq 1}\sum_{k+l=n,~k>0}\frac{\hbar^{-l-k}}{(k-1)!} \sum_{l_1+\cdots + l_r=k,~ l_i>0} \binom{k-1}{l_1,\ldots,l_i - 1,\ldots,l_r}  \langle T_a \psi^l\prod_{j=1}^{r}T_{j}^{l_j}  , T_b\rangle_{0,2,\beta} \nonumber \\
        &\cdot (t_{1}^{l_1} \cdots t_{i-1}^{l_{i-1}} t_{i}^{l_i - 1} t_{i+1}^{l_{i+1}} \cdots t_{r}^{l_r})
        \nonumber
    \end{align}
    We rearrange the first and second terms in \eqref{pro_qde2} without the factor $\sum_{\beta \neq 0}q^\beta \sum_{b}T^b$: 1) For the first term in \eqref{pro_qde2}, we break this into three pieces in terms of $n$ and $l$: i) $n=0$, ii) $n\geq 1,~l=0$, iii) $n\geq 1,~l\geq 1$. For the last piece $n\geq 1,~l \geq 1$, we substitute $l$ by $l+1$. The result is the following:
    \begin{align}
    \label{pro_qde3}
        &\big(\int_{\beta}T_i \big)  \langle T_a  , T_b\rangle_{0,2,\beta} \\
        &+ \big(\int_{\beta}T_i \big) \sum_{n\geq 1} \bigg[
        \frac{\hbar^{-n}}{n!} \sum_{l_1+\cdots + l_r=n} \binom{n}{l_1,\ldots,l_r} \prod_{j=1}^{r}t_{j}^{l_j} \langle T_a \prod_{j=1}^{r}T_{j}^{l_j}  , T_b\rangle_{0,2,\beta}  \nonumber \\
        &+\sum_{k+l=n-1,l\geq 0}\frac{\hbar^{-(l+1)-k}}{k!} \sum_{l_1+\cdots + l_r=k} \binom{k}{l_1,\ldots,l_r}  \prod_{j=1}^{r}t_{j}^{l_j} \langle T_a \psi^{l+1}\prod_{j=1}^{r}T_{j}^{l_j}  , T_b\rangle_{0,2,\beta} 
        \bigg] \nonumber
    \end{align}
    2) For the second term in \eqref{pro_qde2}, we substitute $k$ by $k+1$ and $l_i$ by $l_{i}+1$, respectively. Then,
    \begin{equation}
    \label{pro_qde4}
        \sum_{n\geq 1}\sum_{k+l=n-1,k\geq 0}\frac{\hbar^{-l-(k+1)}}{k!} \sum_{l_1+\cdots + l_r=k} \binom{k}{l_1,\ldots,l_r}  \prod_{j=1}^{r}t_{j}^{l_j}\langle T_a \psi^l (T_{1}^{l_1} \cdots T_{i-1}^{l_{i-1}} T_{i}^{l_i + 1} T_{i+1}^{l_{i+1}} \cdots T_{r}^{l_r})  , T_b\rangle_{0,2,\beta} 
    \end{equation}
    Combining \eqref{pro_qde3} and \eqref{pro_qde4}, one can write \eqref{pro_qde2} as follows
    \begin{align}
        \label{pro_qde5}
        &\sum_{\beta \neq 0} q^\beta \sum_b T^b \bigg[ \big(\int_{\beta}T_i \big)  \langle T_a  , T_b\rangle_{0,2,\beta} \\
        &+ \sum_{n\geq 1} \bigg[
        \frac{\hbar^{-n}}{n!} \sum_{l_1+\cdots + l_r=n} \binom{n}{l_1,\ldots,l_r}  \big(\int_{\beta}T_i \big)\prod_{j=1}^{r}t_{j}^{l_j}
        \langle T_a \prod_{j=1}^{r}T_{j}^{l_j}  , T_b\rangle_{0,2,\beta}  \nonumber \\
        &+\sum_{k+l=n-1}\frac{\hbar^{-l-k-1}}{k!} \sum_{l_1+\cdots + l_r=k} \binom{k}{l_1,\ldots,l_r} \prod_{j=1}^{r}t_{j}^{l_j}
         \nonumber\\
        &\bigg[  \big(\int_{\beta}T_i \big)\langle T_a \psi^{l+1}\prod_{j=1}^{r}T_{j}^{l_j}  , T_b\rangle_{0,2,\beta}
        + \langle T_a \psi^l (T_{1}^{l_1} \cdots T_{i-1}^{l_{i-1}} T_{i}^{l_i + 1} T_{i+1}^{l_{i+1}} \cdots T_{r}^{l_r})  , T_b\rangle_{0,2,\beta}  \bigg]\bigg]\bigg] \nonumber
    \end{align}
    Apply the Divisor equation \eqref{diveqn} to \eqref{pro_qde5}. Then,
    \begin{align}
        \label{pro_qde6}
        &\sum_{\beta \neq 0} q^\beta \sum_b T^b \bigg[  \langle T_a  , T_b \mid \Tilde{T_i} \rangle_{0,2|1,\beta} \\
        &+ \sum_{n\geq 1} \bigg[
        \frac{\hbar^{-n}}{n!} \sum_{l_1+\cdots + l_r=n} \binom{n}{l_1,\ldots,l_r} 
        \prod_{j=1}^{r}t_{j}^{l_j} \langle T_a \prod_{j=1}^{r}T_{j}^{l_j}  , T_b \mid \Tilde{T_i} \rangle_{0,2|1,\beta}  \nonumber \\
        &+\sum_{k+l=n-1}\frac{\hbar^{-l-k-1}}{k!} \sum_{l_1+\cdots + l_r=k} \binom{k}{l_1,\ldots,l_r} \prod_{j=1}^{r}t_{j}^{l_j}
        \langle T_a \psi^{l+1}\prod_{j=1}^{r}T_{j}^{l_j}  , T_b \mid \Tilde{T_i} \rangle_{0,2|1,\beta}
        \bigg]\bigg] \nonumber
    \end{align}

    We compute the right-hand side of \eqref{eqn_qde}: $\Tilde{T}_i \star s_a$. Applying the definition of $H^*(T)$-action \eqref{def_staroper} to the first term of $s_a$ gives us
    \begin{equation}
        \label{pro_qde7}
        \Tilde{T}_i \star (e^{\delta / \hbar} \cup T_a) =  e^{\delta / \hbar} \cup T_a \cup T_i +  \sum_{\beta\neq 0,b}q^\beta T^b \langle e^{\delta / \hbar} \cup T_a, T_b \mid \Tilde{T}_i \rangle .
    \end{equation}
    The first term of \eqref{pro_qde7} is the same as \eqref{pro_qde1}, since the $2|1$-operation recovers the product of the ordinary cohomology at the constant term. Again, apply $\Tilde{T}_i\star$ on the second term of $s_a$ in \eqref{pro_qde_secexpn}. Then, one can simplify as follows
    \begin{align}
        \label{pro_qde8}
        &\Tilde{T}_i\star \bigg( \sum_{\beta'}q^{\beta'} \sum_{\alpha'} T^{\alpha'} \hbar^{-1} \sum_{n\geq 0} \sum_{k+l=n} \frac{\hbar^{-l-k}}{k!} \sum_{l_1+\cdots + l_r=k} \binom{k}{l_1,\ldots,l_r} \prod_{j=1}^{r}t_{j}^{l_j}
        \langle T_a \psi^l\prod_{j=1}^{r}T_{j}^{l_j}  , T_{\alpha'} \rangle_{0,2,\beta'}  \bigg) \\
        &=  \sum_{\beta'}q^{\beta'} \sum_{\alpha'} \sum_{\beta'',b} q^{\beta''}T^b \langle T^{\alpha'}, T_b \mid \Tilde{T}_i \rangle_{0,2|1,\beta''} \nonumber \\
        & \cdot  \sum_{n\geq 0} \sum_{k+l=n} \frac{\hbar^{-l-k-1}}{k!} \sum_{l_1+\cdots + l_r=k} \binom{k}{l_1,\ldots,l_r} \prod_{j=1}^{r}t_{j}^{l_j} \langle T_a \psi^l\prod_{j=1}^{r}T_{j}^{l_j}  , T_{\alpha'} \rangle_{0,2,\beta'} 
        \nonumber \\
        &=\sum_{\beta \neq 0} q^\beta \sum_{\substack{\beta'+\beta''=\beta ,\nonumber \\
        \beta' \neq 0 }}\sum_{b} T^b \sum_{n\geq 0} \sum_{k+l=n} \frac{\hbar^{-l-k-1}}{k!} \nonumber \\
        & \cdot \sum_{l_1+\cdots + l_r=k} \binom{k}{l_1,\ldots,l_r} \prod_{j=1}^{r}t_{j}^{l_j} 
        \sum_{\alpha'}\langle T^{\alpha'}, T_b \mid \Tilde{T}_i \rangle_{0,2|1,\beta''}  \langle T_a \psi^l\prod_{j=1}^{r}T_{j}^{l_j}  , T_{\alpha'} \rangle_{0,2,\beta'} \nonumber
    \end{align}
    Applying the topological recursive relation \eqref{trr}, adding \eqref{pro_qde6} and the second term of \eqref{pro_qde7} are identical to \eqref{pro_qde8}. Therefore, we verified that the quantum differential equation holds.
\end{proof}

The idea of using the dual connection below was introduced from \cite[\S2.3]{quamtumhyper}, due to Bumsig Kim.
\begin{Def}
With the same notation used to define the quasimap Givental connection, its \textit{dual connection} $\Tilde{\nabla}$ is defined by
\[
    \Tilde{\nabla}_{\frac{\partial}{\partial t_i}} \big(\sum_{j=0}^{m}a_jT_j\big):=\hbar \sum_{j=0}^{m}\frac{\partial a_j}{\partial t_i}T_j + \sum_{j=0}^{m}a_j\Tilde{T_i} \star T_j.
\]
\end{Def}
For simplicity, denote the connection $\Tilde{\nabla}_{\frac{\partial}{\partial t_i}}$ by $\Tilde{\nabla}_i$. Write $q_i=e^{t_i}$ and $q^{\beta} = e^{\sum_{j=1}^{r} t_j\int_{\beta}T_j} $. We prove a Leibniz rule for the operator $\hbar {\partial } / {\partial t_i}$.
\begin{Lem}
\label{pro_leibniz}
    For $H^*(X_\Sigma) \otimes_{\C} \C[[\text{Eff}]] $-valued functions $G$ and $H$ in $t_0,\ldots,t_r$, the following equation holds
    \begin{equation*}
        \hbar \frac{\partial }{\partial t_i}\langle G, H\rangle  = \langle \nabla_i G , H \rangle + \langle G, \Tilde{\nabla}_i H \rangle, 
    \end{equation*}
    where $i= 0,\ldots, r$ and the $\langle -, -\rangle$ the cup product.
\end{Lem}
\begin{proof}
    Write
    \begin{align*}
        G=\sum_{\beta \in \text{Eff}} \sum_{i=0}^{r} z_{i,\beta} T_i e^{\sum_{k=0}^{r} t_k \int_{\beta}T_k}, \quad
        H=\sum_{\beta' \in \text{Eff}} \sum_{j=0}^{r} w_{j,\beta'} T_j e^{\sum_{l=0}^{r} t_l \int_{\beta'}T_l},
    \end{align*}
    where $z_{i,\beta},w_{j,\beta'} \in \C$. Then, one can see that
    \begin{align}
        \label{eqn_lhsderiv}
        \frac{\partial}{\partial t_p } \langle G, H \rangle &= \frac{\partial}{\partial t_p } \sum_{i,j} T_i \cup T_j \sum_{\beta ,\beta'}  z_{i,\beta} w_{j,\beta'} e^{\sum_{k=0}^{r} t_k \int_{\beta + \beta'}T_k } \nonumber \\
        &=\sum_{i,j} T_i \cup T_j \sum_{\beta ,\beta'}  z_{i,\beta} w_{j,\beta'} \big( \int_{\beta + \beta'}T_p \big) e^{\sum_{k=0}^{r} t_k \int_{\beta + \beta'}T_k }
    \end{align}

    On the other hand, by definition, one can simplify $\nabla_p G$
    \begin{align}
        \label{eqn_nablaG}
        \nabla_p G &= \nabla_p \bigg( \sum_{\beta } \sum_{i=0}^{r} z_{i,\beta} T_i e^{\sum_{k=0}^{r} t_k \int_{\beta}T_k} \bigg) \nonumber \\
        &=\hbar\sum_{i=0}^{r} \frac{\partial}{\partial t_p }\sum_{\beta }  z_{i,\beta} T_i e^{\sum_{k=0}^{r} t_k \int_{\beta}T_k} - \sum_{i=0}^{r}\bigg( \sum_{\beta }  z_{i,\beta}  e^{\sum_{k=0}^{r} t_k \int_{\beta}T_k}  \bigg) \Tilde{T}_p \star T_i \nonumber \\
        &= \hbar \sum_{i=0}^{r} \sum_{\beta }  z_{i,\beta} T_i \big(\int_{\beta}T_p \big) e^{\sum_{k=0}^{r} t_k \int_{\beta}T_k} \nonumber \\
        - \sum_{i=0}^{r}&\bigg( \sum_{\beta }  z_{i,\beta}  e^{\sum_{k=0}^{r} t_k \int_{\beta}T_k}  \bigg) \bigg( \sum_{\beta' } \sum_{j=0}^{r} e^{\sum_{k'=0}^{r} t_{k'} \int_{\beta'}T_{k'}} \langle T_i,T_j | \Tilde{T}_p  \rangle_{0,2|1,\beta'} T^j  \bigg) \nonumber  \\
        &= \hbar \sum_{i=0}^{r} \sum_{\beta }  z_{i,\beta} T_i \big(\int_{\beta}T_p \big) e^{\sum_{k=0}^{r} t_k \int_{\beta}T_k} \nonumber \\
        - \sum_{i=0}^{r}& \sum_{j=0}^{r} \sum_{\beta, \beta'}   T^j z_{i,\beta} e^{\sum_{k=0}^{r} t_k \int_{\beta + \beta'}T_k} \langle T_i,T_j | \Tilde{T}_p  \rangle_{0,2|1,\beta'}   .
    \end{align}
    Similarly, we have
    \begin{align}
    \label{eqn_nablaH}
     \Tilde{\nabla}_p H  &= \hbar \sum_{j=0}^{r} \sum_{\beta' }  w_{j,\beta'} T_j \big(\int_{\beta'}T_p \big) e^{\sum_{k=0}^{r} t_k \int_{\beta'}T_k} \nonumber \\
        + \sum_{j=0}^{r}& \sum_{l=0}^{r} \sum_{\beta, \beta'}   T^l w_{j,\beta'} e^{\sum_{k=0}^{r} t_k \int_{\beta + \beta'}T_k} \langle T_j,T_l | \Tilde{T}_p  \rangle_{0,2|1,\beta}   .
    \end{align}
    Thus, from \eqref{eqn_nablaG} and \eqref{eqn_nablaH}, a straightforward calculation and reindexing subscripts show that
    \begin{align}
    \label{eqn_rhs_twonabla}
        &\langle \nabla_p G, H \rangle + \langle  G, \Tilde{\nabla}_p H \rangle \nonumber  \\
        =  \hbar &\sum_{i=0}^{r} \sum_{j=0}^{r} \sum_{\beta,\beta' }  z_{i,\beta}  w_{j,\beta'} T_i \cup T_j \big(\int_{\beta}T_p \big) e^{\sum_{k=0}^{r} t_k \int_{\beta + \beta'}T_k}  \nonumber \\
        - \sum_{i=0}^{r}& \sum_{j=0}^{r} \sum_{k=0}^{r} \sum_{\beta, \beta', \beta''}   T^j \cup T^k  z_{i,\beta} w_{k,\beta''} e^{\sum_{k=0}^{r} t_k \int_{\beta + \beta' + \beta''}T_k} \langle T_i,T_j | \Tilde{T}_p  \rangle_{0,2|1,\beta'} \nonumber  \\
        +  \hbar &\sum_{i=0}^{r} \sum_{j=0}^{r} \sum_{\beta,\beta' }  z_{i,\beta}  w_{j,\beta'} T_i \cup T_j \big(\int_{\beta'}T_p \big) e^{\sum_{k=0}^{r} t_k \int_{\beta + \beta'}T_k}  \nonumber \\
        - \sum_{i=0}^{r}& \sum_{j=0}^{r} \sum_{k=0}^{r} \sum_{\beta, \beta', \beta''}   T^j \cup T^k  z_{i,\beta} w_{k,\beta''} e^{\sum_{k=0}^{r} t_k \int_{\beta + \beta' + \beta''}T_k} \langle T_i,T_j | \Tilde{T}_p  \rangle_{0,2|1,\beta'} \nonumber \\
        =  \hbar &\sum_{i=0}^{r} \sum_{j=0}^{r} \sum_{\beta,\beta' }  z_{i,\beta}  w_{j,\beta'} T_i \cup T_j \big(\int_{\beta+\beta'}T_p \big) e^{\sum_{k=0}^{r} t_k \int_{\beta + \beta'}T_k} .
    \end{align}
    Since \eqref{eqn_lhsderiv} and \eqref{eqn_rhs_twonabla} are the same, we verified the statement.
\end{proof}

Observe how the dual connection acts on the fundamental class in $M=H^0(X_\Sigma) \oplus H^2(X_\Sigma)$.
\begin{Lem}
    \label{pro_dualconnone}
    The following equation holds:
    \[
    \Tilde{\nabla}_{i_k} \cdots \Tilde{\nabla}_{i_1} 1 = \Tilde{T}_{i_k}\star \cdots \star \Tilde{T}_{i_1} \star 1 + \hbar A(T,q,\hbar),
    \]
    for some formal power series $A(T,q,\hbar)$ in $T_i,~q_j$, and $\hbar$.
\end{Lem}
\begin{proof}
    We prove by using induction. Suppose $k=1$. Without loss of generality, we put $i:=i_1$. By the definition of the dual connection, one can write
    \[
    \Tilde{\nabla}_{i} 1=\Tilde{T}_i\star 1.
    \]
    
    We proceed the inductive case. Using the equation of the module structure \eqref{eqn_WDVV} and the formula $\frac{\partial}{\partial t_i }q^{\beta} = q^{\beta}\int_{\beta}T_i $, one can compute
    \begin{align}
    \label{eqn_dualconn1}
         & \Tilde{\nabla}_{i_k} \big( \Tilde{T}_{i_{k-1}} \star \cdots \star \Tilde{T}_{i_{1}} \star 1 \big) \\
         =&\hbar \sum_{j} \frac{\partial}{\partial t_{i_k}} \sum_{\beta}q^\beta T^{j} \langle 1 , T_l  \mid \Tilde{T}_{i_{k-1}} \cup \cdots \cup \Tilde{T}_{i_{1}} \rangle_{0,2|1,\beta} + \Tilde{T}_{i_{k}} \star \cdots \star \Tilde{T}_{i_{1}} \star 1 \nonumber \\
         =& \hbar \sum_{j} \bigg( \int_{\beta}T_{i_k} \bigg)  \sum_{\beta}q^\beta T^{j} \langle 1 , T_l  \mid \Tilde{T}_{i_{k-1}} \cup \cdots \cup \Tilde{T}_{i_{1}} \rangle_{0,2|1,\beta} + \Tilde{T}_{i_{k}} \star \cdots \star \Tilde{T}_{i_{1}} \star 1. \nonumber
    \end{align}
    On the other hand, for a formal power series
    \[A(T,q,\hbar) = \sum_{i} T_i \sum_{\beta}  \sum_{n}  z_{i n \beta} \hbar^n q^\beta, \]
    similar computation shows that
    \begin{align}
    \label{eqn_dualconn2}
         & \Tilde{\nabla}_{i_k} A(T,q,\hbar) \\
         =& \hbar \sum_{i} T_i \sum_{\beta}  \sum_{n}  z_{i n \beta} \hbar^n q^\beta  \bigg( \int_{\beta}T_{i_k} \bigg) + \sum_{i} \Tilde{T}_{i_k} \star T_i \sum_{\beta}  \sum_{n}  z_{i n \beta} \hbar^n q^\beta \nonumber\\
         =& \hbar  \sum_{i} T_i \sum_{\beta}  \sum_{n}  z_{i n \beta} \hbar^n q^\beta  \bigg( \int_{\beta}T_{i_k} \bigg) + \sum_{i} \sum_{l,\beta'}q^{\beta'} T^l  \langle T_i, T_l  \mid \Tilde{T}_{i_k} \rangle_{0,2|1,\beta'}  \sum_{\beta}  \sum_{n}  z_{i n \beta} \hbar^n q^\beta .\nonumber
    \end{align}
    From \eqref{eqn_dualconn1} and \eqref{eqn_dualconn2}, we complete our proof of the statement.
\end{proof}

\section{$I$-function and relations in the Quantum module}
\label{sec_Iftn}
We define a cohomological-valued function via quasimap invariants by following \cite[\S10.3]{coxmirror}.
\begin{Def}
Define the \textbf{$I$-function} in the following way
\begin{equation}
    \label{eqn_IftnDef}
    I:=\sum_{a=0}^{m} T^a \langle s_a, 1 \rangle,
\end{equation}
where $\langle \alpha, \beta \rangle := \int_{X_\Sigma} \alpha \cup \beta$.
\end{Def}
The $I$-function has the following expression by a straightforward computation
\begin{equation}
        \label{eqn_IftnExp}
        I= e^{(t_0 + \delta)/\hbar}\bigg( 1 + \sum_{\beta \neq 0}\sum_{a=1}^{m}q^\beta \langle \frac{T_a}{\hbar-\psi}, 1 \rangle_{0,2,\beta}T^a \bigg) .
\end{equation}

We prove the following theorem, which gives us a relation in the quantum module of $X_\Sigma$.
\begin{Thm}
    \label{pro_relation}
    Let $P(\hbar\partial/\partial t, e^{t}, \hbar)$ be a formal power series in
    \[\hbar \frac{\partial}{\partial t_0}, \ldots, \hbar \frac{\partial}{\partial t_r}, e^{t_0}, \ldots, e^{t_r}, \hbar, \]
    and $P(\Tilde{T}, q, 0)$ be the formal power series in $H^*([V/T]) \otimes_{\C} \C[[\text{Eff}]] $ obtained by replacing
    \[
    \hbar \frac{\partial}{\partial t_i} \rightarrow \Tilde{T}_i,~e^{t_j} \rightarrow q_j,~ \hbar \rightarrow 0,
    \]
    and the composition of differential operators by the action in \eqref{def_staroper}. If \[P(\hbar\partial/\partial t, e^{t}, \hbar)I:= \sum_a T^a P(\hbar\partial/\partial t, e^{t}, \hbar)\langle s_a, 1 \rangle =0,\] then the relation \[P(\Tilde{T}, q, 0)\star 1=0\] holds in the quantum module $QM(X_\Sigma)$.
\end{Thm}

\begin{proof}
Since $\{T_a\}$ forms a basis of $H^*(X_\Sigma)$, the hypothesis $PI=0$ implies 
\begin{equation}
\label{eqn_Pseczero}
P\langle s_a, 1 \rangle =0        
\end{equation}
for all $a=0,\ldots, m$. With the result of Proposition \ref{pro_qde}, repetitive uses of Lemma \ref{pro_leibniz} and Lemma \ref{pro_dualconnone} to \eqref{eqn_Pseczero} give us
\[
\langle s_a, P(\Tilde{T},q,0) \star 1 \rangle + \hbar \langle s_a, A(T,q,\hbar)\rangle=0,
\]
for some formal power series $A(T,q,\hbar)$. Since $\hbar$ is a free variable and the fact that $s_a$ form a basis, we conclude that there is a relation given by $P(\Tilde{T},q,0) \star 1$=0 in the quantum module $QM(X_\Sigma)$.
\end{proof}

\section{The Batyrev modules for semipositive spaces}
\label{sec_bat}
\subsection{Batyrev modules} We introduce the Batyrev ring for a smooth projective toric variety $X_\Sigma$ defined in \cite[Def 5.1]{batyrev}. Define the following two ideals in $\C[x_\rho | \rho \in \Sigma(1)] \otimes_{\C} \C[[\text{NE}_\Z]]$
\begin{align*}
    P_\Sigma&:=\big\langle \sum_{\rho \in \Sigma(1)}\langle m,u_\rho \rangle x_\rho :m\in M \big\rangle,\\
    SR_\Sigma&:= \big\langle x_{\rho_{i_1}}\cdots x_{\rho_{i_k}} - q^{\beta_P}x_{\rho_{k_1}}\cdots x_{\rho_{j_l}} : P=\{ \rho_{i_1},\cdots ,\rho_{i_k} \}\text{ is a primitive collection}    \big\rangle,
\end{align*}
where the generators of the Stanley-Reisner ideal $SR_\Sigma$ come from the primitive relation \eqref{eqn_primitiverelation} of $P$ and the primitive class $\beta_P$.The \textbf{Batyrev ring} of $X_\Sigma$ is defined by
    \[QH_{\text{Bat}}^*(X_\Sigma):= \C[x_\rho|\rho \in \Sigma(1)]  \otimes_{\C}  \C[[\text{NE}_\Z]] / (P_\Sigma + SR_\Sigma ) . \]
The Batyrev ring $QH_{\text{Bat}}^*(X_\Sigma)$ has a natural $\C[x_\rho|\rho \in \Sigma(1)]  \otimes_{\C}  \C[[\text{NE}_\Z]]$-module structure given by the quotient map. Denote this module by Bat$M(X_\Sigma)$ and call the \textbf{Batyrev module} of $X_\Sigma$.

The Mori cone and the quasimap-effective cone can be identified.
\begin{Pro}
    \label{pro_morieff}
    For a smooth projective toric variety $X_\Sigma$,
    \begin{equation*}
        \text{Eff}  = \text{NE}_\Z
    \end{equation*}
\end{Pro}
\begin{proof}
    ($\subseteq$) This part is from \cite[Lem 3.1.3]{toricquasi}

    ($\supseteq$) From the description of the Mori cone in terms of the extremal rays in \eqref{eqn_moriprimitive}, it is enough to show that for a primitive collection $P$, the primitive class $\beta_P$ is quasimap-effective. One can write $\beta_P =(b_\rho)_{\rho\in \Sigma(1)} \in \R^{\Sigma(1)}$ from the exact sequence \eqref{eqn_toricses2} (or see \cite[Def 6.4.10]{coxtoric}). Let $\gamma$ be the minimal cone containing $\sum_{\rho \in P}u_\rho$. We define a quasimap in $Q_{0,m|k}(X_\Sigma, \beta)$ with $m\geq 2$ whose source curve is $\PP^1$.
    Pick the line bundles $L_\rho =\calO_{\PP^1}(b_\rho)$ on $\PP^1$. Then, $\sum_\rho b_\rho u_\rho =0$ from the exact sequence \eqref{eqn_toricses2}. Thus, the trivialization condition in Remark \ref{remark}(1) is satisfied. Since $\beta_P$ is effective, the stability condition of a quasimap holds. From \eqref{eqn_primitivecoeff}, one can choose sections as follows
    \begin{equation}
        \label{eqn_secprimitive}
        s_\rho = \begin{cases}
            \ell &\rho \in P \\
            0 &\rho \in \gamma(1) \\
            z_\rho &\text{otherwise},
        \end{cases}
    \end{equation}
    where $\ell$ is a homogeneous polynomials of degree $1$ and $z_\rho$ are nonzero complex numbers. One can see that the sections $(s_\rho)_{\rho \in \Sigma(1)}$ can have only finitely many base points including the root $\mathbf{x}$ of $\ell$. Write $B$ as the finite set of base points including $\mathbf{x}$. To meet the nondegeneracy condition of a quasimap, since the fan $\Sigma$ is complete, there exist a maximal cone $\sigma_{\text{max}}$ containing $\gamma$. From the definition of $(s_\rho)_{\rho \in \Sigma(1)}$, one can see that for arbitrary $\mathbf{y} \in \PP^1 \backslash B$, $s_\rho(\mathbf{y}) \neq 0$ for all $\rho \not\subset \sigma_{\text{max}}$. Therefore, our choice of line bundles and sections define a quasimap to $X_\Sigma$.
    
\end{proof}

\section{Proof of Theorem \ref{conj}}
\label{sec_conj}
We prove Conjecture \ref{conj} by relating three cohomological-valued functions. The first function is the $I$-function defined in \eqref{eqn_IftnDef}.
\subsection{$I^{\text{GKZ}}$-functions} The second function arises from the GKZ system of $X_\Sigma$ as follows, see \cite[(11.86) in p389]{coxmirror},
\begin{equation}
    \label{eqn_IftnDef2}
    I^{\text{GKZ}}:=e^{(t_0 + \delta)/\hbar} \sum_{\beta \in \text{NE}_{\Z}(X_\Sigma) } q^\beta \prod_{\rho \in \Sigma(1)} \frac{\prod_{m=-\infty}^{0} (D_{\rho} + m\hbar ) }{\prod_{m=-\infty}^{ \int_{\beta} D_{\rho}  }  (D_{\rho} + m\hbar ) }.
\end{equation}
We call it the \textit{GKZ $I$-function}. Recall that a smooth projective toric variety is \textit{semipositive} if its anti-canonical divisor is numerically effective, i.e., $-K_{X_\Sigma}\cdot [C] \geq 0$ for all effective curve classes $[C]$. We compute the constant term in the expansion of $I^{\text{GKZ}}$ with respect to $1/\hbar$.
\begin{Pro}
    \label{pro_coeff}
    The GKZ $I$-function for a smooth projective toric semipositive variety $X_\Sigma$ has the following expansion
    \begin{equation}
    \label{eqn_GZKexpand}
        I^{\text{GKZ}}=e^{(t_0 + \delta)/\hbar} (I_0+ I_1\hbar^{-1} + o(\hbar^{-2}) )
    \end{equation}
    with $I_0=1$ and for some $I_1$.
\end{Pro}
\begin{proof}
    Fix an effective curve class $\beta \in H_2(X_\Sigma, \Z)$ and divide into four cases with respect to $\rho$. In each case, we determine the highest power of $\hbar^{-1}$ in the expansion of \[\Phi(\rho) := \frac{\prod_{m=-\infty}^{0} (D_{\rho} + m\hbar ) }{\prod_{m=-\infty}^{ \int_{\beta} D_{\rho}  }  (D_{\rho} + m\hbar ) },\]
    with respect to $\hbar^{-1}$.

    (case 1: $D_\rho \cdot \beta =0$) There is no contribution because of cancellation. We can write $(\hbar^{-1})^{D_\rho \cdot \beta} = 1$ for the highest power.

    (case 2: $D_\rho \cdot \beta =-1$) In this case, $\Phi(\rho)=D_\rho$ and $(\hbar^{-1})^{D_\rho \cdot \beta +1} = 1$.

    (case 3: $D_\rho \cdot \beta \leq -2$) One can compute
    \[\Phi(\rho) = D_\rho ((\hbar^{-1})^{D_\rho \cdot \beta +1}(-1)^{-(D_\rho \cdot \beta +1)}(-D_\rho\cdot\beta + 1) + \cdots + D_\rho^{-D_\rho\beta -1}   ). \]
    Thus, the highest power is $D_\rho\cdot \beta + 1$.
    
    (case 4: $D_\rho \cdot \beta > 0$) The expansion of $\Phi(\rho)$ with respect to $\hbar^{-1}$ is
    \[\Phi(\rho) = \frac{1}{(D_\rho \cdot \beta)!} (\hbar^{-1})^{D_\rho \cdot \beta}\prod_{m=1}^{D_\rho \cdot \beta}\sum_{l \geq 0}\big(-\frac{D_\rho}{m\hbar} \big)^{l}. \]
    The highest power is $D_{\rho}\cdot \beta$.

    Multiplying all the highest powers gives us
    \[N:=\sum_{\rho \in \Sigma(1)}D_\rho \cdot \beta + |\{\rho \mid D_\rho \cdot \beta \leq -1\}|= -K_{X_\Sigma} \cdot \beta + |\{\rho \mid D_\rho \cdot \beta \leq -1\}| . \]
    Since $X_\Sigma$ is semipositive, $N \geq 0$. To determine $I_0$, we consider the case $N=0$. In this case, both terms in $N$ are zero. The latter term being equal to zero implies that $D_\rho \cdot \beta \geq 0 $ for all $\rho \in \Sigma(1)$. Then, we must have $D_\rho \cdot \beta =0$ for all $\rho \in \Sigma(1)$, because $-K_{X_\Sigma}\cdot \beta =0$. Recall that the pairing between $H_2(X_\Sigma)$ and $H^2(X_\Sigma)$ is perfect. Since $D_\rho$ generate $H^2(X_\Sigma)$, $\beta$ must be zero. Hence, we conclude that $I_0=1$.
\end{proof}

\subsection{Restricted big $J$-function} The third function we introduce is the big $J^{\epsilon}$-function $J^{\epsilon}(\mathbf{t},z)$ from \cite[Def 5.1.1]{genus0}. We only consider the case $\epsilon=0+$ and where the generic insertions $\mathbf{t} \in H^*(X_\Sigma)$ are restricted to $\mathbf{0}$. Then,
\begin{equation*}
    J^{0+}\lvert_{\textbf{t}=0}=1+\sum_{\beta \neq 0} q^{\beta}(\text{ev}_{\bullet})_{*}\frac{ [Q_{0,0+\bullet}(X_\Sigma,\beta)_0]^{\text{vir} } }{e^{\C^*}(N_{F_0/QG})},
\end{equation*}
where $QG$ is notation for a quasimap graph spaces with $\C^*$-action, $N_{F_0/QG}$ is the virtual normal bundle of a fixed locus $F_0 \simeq Q_{0,0+\bullet}(X_\Sigma,\beta)_0$ whose $\PP^1$-parametrized component has a degree $\beta$ base point at 0. We refer readers to \cite{genus0} for precise definitions. In \cite{toricquasi}, they compute the equivariant eular characteristic $e^{\C^*}(N_{F_0/QG})$
\[
\frac{1}{e^{\C^*}(N_{F_0/QG})}=\prod_{\rho \in \Sigma(1)} \frac{\prod_{m=-\infty}^{0} (D_{\rho} + m\hbar ) }{\prod_{m=-\infty}^{ \int_{\beta} D_{\rho}  }  (D_{\rho} + m\hbar ) }.
\]
Therefore, by identifying the Mori cone and the quasimap effective cone through Proposition \ref{pro_morieff}, one can see
\begin{equation}
    \label{eqn_gkzJftn}
    I^{\text{GKZ}} = e^{(t_0 + \delta)/\hbar}J^{0+}\lvert_{\mathbf{t}=0}.
\end{equation}

We write
\[
J^{0+}\lvert_{\mathbf{t}=0} = J_{0}^{0+}(q) + \frac{1}{\hbar}J_{1}^{0+} + O(\frac{1}{\hbar^2}).
\]
Proposition \ref{pro_coeff} gives us
\begin{equation}
\label{eqn_Jis1}
    J_{0}^{0+} = 1.
\end{equation}
Combining \eqref{eqn_Jis1} with Corollary 5.5.1 in \cite{genus0}, we identify the $I$-function expressed in \eqref{eqn_IftnExp} with \eqref{eqn_gkzJftn}.
\begin{Pro}
The three cohomological-valued functions $I^{GKZ}$, $e^{(t_0 + \delta)/\hbar}J^{0+}\lvert_{\mathbf{t}=0}$ and $I$ are identical
    \label{IandJ}
    \begin{equation*}
    I= e^{(t_0 + \delta)/\hbar}J^{0+}\lvert_{\mathbf{t}=0} = I^{\text{GKZ}}.
\end{equation*}
\end{Pro}

\subsection{Isomorphism between BatM($X_\Sigma$) and $QM(X_\Sigma)$}

For $\beta \in H_2(X_\Sigma,\Z) $, define an operator 
\begin{align}
    \label{eqn_GKZoper}
    \square_{\beta,\hbar}&:=\prod_{\rho:\int_\beta D_\rho > 0} \prod_{m=0}^{\int_\beta D_\rho-1} \bigg( \bigg( \sum_{j} \int_{\beta_j}D_\rho \hbar\frac{\partial}{\partial t_j} \bigg) - m\hbar  \bigg)   \\
    &- q^\beta \prod_{\rho:\int_\beta D_\rho < 0} \prod_{m=0}^{-\int_\beta D_\rho-1} \bigg( \bigg( \sum_{j} \int_{\beta_j}D_\rho  \hbar\frac{\partial}{\partial t_j} \bigg) - m\hbar  \bigg). \nonumber
\end{align}
From \cite[\S11.2 p391]{coxmirror}, the equation in Proposition \ref{IandJ} implies that the quantum differential operator \eqref{eqn_GKZoper} annihilates both the GKZ $I$-funtion $I^{\text{GKZ}}$ and the quasimap $I$-function
\begin{equation*}
    \square_{\beta,\hbar} I^{\text{GKZ}}= \square_{\beta,\hbar} I = 0.
\end{equation*}
Thus, with the unique lift $\Tilde{T}_\rho$ in $H^2(T)$, Theorem \ref{pro_relation} gives us relations in the quantum module $QM(X_\Sigma)$
\begin{equation}
    \label{eqn_relationinQM}
    \bigg( \prod_{\rho:\int_\beta D_\rho > 0} \Tilde{T}_{\rho}^{\int_{\beta}D_\rho} - q^\beta \prod_{\rho:\int_\beta D_\rho < 0} \Tilde{T}_{\rho}^{-\int_{\beta}D_\rho} \bigg) \star 1 = 0  .
\end{equation}

Recall that we can identify $H^0(T) \oplus H^2(T)$ with $H^0(X_\Sigma) \oplus H^2(X_\Sigma)$ via the Kirwan map. Take a homogeneous basis $\{\sigma_0:=1,\sigma_1,\ldots,\sigma_r\}$ of $H^0(T) \oplus H^2(T)$. Define a map $\varphi$ from $\C[[\text{NE}_\Z]][\sigma_1,\ldots, \sigma_r]$ to $QM(X_\Sigma)$ by sending as follows
\[\sigma_i \longmapsto \sigma_i \star 1.\]
Since the relation \eqref{eqn_relationinQM} holds in both BatM$(X_{\Sigma})$ and $QM(X_\Sigma)$, the map $\varphi$ extends to a $\C[[\text{NE}_\Z]][\sigma_1,\ldots, \sigma_r]$-morphism
\[
\varphi: \text{BatM}(X_\Sigma) \longrightarrow QM(X_\Sigma).
\]
Take a basis for $H^*(X_\Sigma)$ including the images of $\sigma_i$. Then, we can represent $\phi$ as a matrix whose entries are in $\C[[\text{NE}_\Z]]$. Since the $H^*(T)$-action \eqref{def_staroper} recovers the ordinary intersection product in the degree zero part, one can show that
\[
\det(\varphi) = 1 + m,
\]
where $m$ is an element in the unique maximal ideal of the local ring $\C[[\text{NE}_\Z]]$. Thus, since $\C[[\text{NE}_\Z]]$ is a power series ring, $\varphi$ is an invertible element. Hence, we proved Conjecture \ref{conj} by showing that there is an isomorphism
\[\text{BatM}(X_\Sigma) \simeq QM(X_\Sigma),\]
for a smooth projective toric semipositive variety $X_\Sigma$.
{
\bibliographystyle{plain}
\bibliography{References.bib}
}
\end{document}